\documentclass{llncs}
\begin{document}
\pagenumbering{arabic}
\pagestyle{headings}
\def\sof{\hfill\rule{2mm}{2mm}}

\title{{\sc RESTRICTED SINGLE OR DOUBLE SIGNED PATTERNS}}
   
\author{{\sc Toufik Mansour}}
\institute{Department of Mathematics,\\ 
	   University of Haifa, Israel 31905\\
	   tmansur@study.haifa.ac.il}

\maketitle
\section*{\centering{\sc Abstract}}

Let $E_n^r=\{[\tau]_a=(\tau_1^{(a_1)},\dots,\tau_n^{(a_n)})| \tau\in S_n,\ 1\leq a_i\leq r\}$ 
be the set of all signed permutations on the symbols $1,2,\dots,n$ with signs $1,2,\dots,r$. 
We prove, for every $2$-letter signed pattern $[\tau]_a$, that the 
number of $[\tau]_a$-avoiding signed permutations in $E_n^r$ is given by 
the formula $\sum\limits_{j=0}^n j!(r-1)^j{n\choose j}^2$. Also we prove that there 
are only one Wilf class for $r=1$, four Wilf classes for $r=2$, and six 
Wilf classes
 for $r\geq 3$.\\
\\
{\bf Key words:} restricted permutations, pattern avoidance, signed permutations.

%===========================================================================
\section*{\centering{\sc 1. Introduction}} 

\ \ \ \ Pattern avoidance proved to be a useful language in a variety of seemingly unrelated 
problems, from stack sorting \cite{Kn,Rt} to the theory of Kazhdan-Lusztig 
polynomials ~\cite{Fb}, and singularities of Schubert varieties \cite{LS,SCb}. 
Signed pattern avoidance proved to be a useful 
language in combinatorial statistics defined in type-$B$ noncrossing partitions, 
enumerative combinatorics, algebraic combinatorics, and geometric combinatorics 
\cite{Rs,MbRs,Cm,Vr}. \\

{\bf Restricted permutations.}
Let $\pi\in S_n$ and $\tau\in S_k$ be two permutations. An {\it occurrence} 
of $\tau$ in $\pi$ is a subsequence $1\leq i_1<i_2<\dots<i_k\leq n$ such that 
$(\pi_{i_1},\dots,\pi_{i_k})$ is order-isomorphic to $\tau$; in such a 
context $\tau$ is usually called a {\it pattern}. We say that $\pi$ {\it avoids} 
$\tau$, or is $\tau$-{\it avoiding}, if there is no occurrence of $\tau$ in $\pi$.
The set of all $\tau$-avoiding permutations in $S_n$ is denoted 
$S_n(\tau)$. 
For an arbitrary finite collection of patterns $T$, we say that $\pi$
avoids $T$ if $\pi$ avoids any $\tau\in T$; the corresponding subset of $S_n$
is denoted $S_n(T)$.\\

{\bf Restricted signed permutations.} 
We say that $(\tau_1^{(a_1)},\dots,\tau_n^{(a_n)})$ is a {\it signed permutation\/} 
and denote it 
by $[\tau]_a$ if $(\tau_1,\dots,\tau_n)\in S_n$ and $a\in [r]^n$. 
In this context, we call $a_1,\dots,a_n$ the 
{\it signs\/} of $\tau$, and we call $\tau_1,\dots,\tau_n$ the 
{\it symbols\/} of $\tau$.

The set of all signed permutations with symbols $a_1,\dots,a_n$ and signs 
$d_1,d_2,\dots,d_r$ we denote by $E_{a_1,\dots,a_n}^{d_1,\dots,d_r}$;  
also we denote $E_n^r=\{ [\tau]_a | \tau\in S_n, 1\leq a_i\leq r \}$. 
Clearly, by definitions $|E_n^r|=n!\cdot r^n$.\\

Similarly to the symmetric group $S_n$ which is generated by the adjacent transpositions $\sigma_i$ for 
$1\leq i\leq n$, where $\sigma_i$ interchanges positions $i$ and $i+1$ (see also the hyperoctahedral group $B_n$ ~\cite{Rs}), the set $E_n^r$ is a  
group which is generated by the adjacent transpositions $\sigma_i$ for $1\leq i\leq n$, 
along with $\sigma_0$ which acts on the right by increasing the first sign; that is, 
	$$(\tau_1^{(a_1)},\tau_2^{(a_2)} ,\dots,\tau_n^{(a_n)})\sigma_0=
	 (\tau_1^{1+(a_1+1(mod\ r))},\tau_2^{(a_2)} ,\dots,\tau_n^{(a_n)}).$$

\begin{example}
The set of all signed permutations with two symbols $1,2$ and two signs $1,2$ is the following set:
$$\begin{array}{ll}
 E_2^2 =\{ &(1^{(1)},2^{(1)}),\ (1^{(1)},2^{(2)}),\ (1^{(2)},2^{(1)}),\ (1^{(2)},2^{(2)}),\\ 
           &(2^{(1)},1^{(1)}),\ (2^{(1)},1^{(2)}),\ (2^{(2)},1^{(1)}),\ (2^{(2)},1^{(2)})\}.
\end{array}$$ 
\end{example}

Let $[\tau]_a\in E_k^r$, and $[\alpha]_b\in E_n^r$; 
we say that $[\alpha]_b$ {\em avoids} $[\tau]_a$ (or is $[\tau]_a$-avoiding) if there is 
no sequence of $k$ indices, $1\leq i_1<i_2<\dots<i_k\leq n$ such that the following two 
conditions hold:

\begin{enumerate}
\item[(i)] 	$(\alpha_{i_1},\dots,\alpha_{i_k})$ is order-isomorphic to $\tau$; 

\item[(ii)] 	$b_{i_j}=a_j$ for all $j=1,2,\dots,k$.
\end{enumerate}

Otherwise, we say that $[\alpha]_b$ {\em contains} $[\tau]_a$ (or is  $[\tau]_a$-containing). 
The set of all $[\tau]_a$-avoiding signed permutations in $E_n^r$ denoted by $E_n^r([\tau]_a)$, 
and in this context $[\tau]_a$ is called a {\it signed pattern\/}. 
For an arbitrary finite collection of signed patterns $T$, we say that $[\alpha]_b$
avoids $T$ if $[\alpha]_b$ avoids any $[\tau]_a\in T$; the corresponding subset of $E_n^r$
is denoted $E_n^r(T)$. 

\begin{example}
As an example, $\Phi=(3,2,1)_{(1,2,2)}=(3^{(1)},2^{(2)},1^{(2)})\in E_3^2$ avoids 
$(2^1,1^1)$; that is, $\Phi\in E_3^2((2^{(1)},1^{(1)}))$. 
\end{example}

Let $T_1$, $T_2$ be two subsets of signed patterns; we say that $T_1$ and $T_2$ are in the same 
{\it $d$-Wilf class} if $|E_n^r(T_1)|=|E_n^r(T_2)|$ for $n\geq 0$, $r\geq d$.\\

In the symmetric group $S_n$, for every $2$-letter pattern $\tau$ the number 
of $\tau$-avoiding permutations is one, and for every pattern $\tau\in S_3$ 
the number of $\tau$-avoiding permutations is given by the Catalan number 
~\cite{Kn}. Also Simion ~\cite{Rs} proved there are similar results for the
hyperoctahedral group $B_n$. Here we are looking for similar results for
$E_n^r$. We show that for every $2$-letter signed pattern $[\tau]_a$ the number of
$[\tau]_a$-avoiding signed permutations in $E_n^r$ is given by 
$\sum\limits_{j=0}^n j!(r-1)^j{n\choose j}^2,$ which generalize the results 
of ~\cite{Rs} (see section $3$).\\

The paper is organized as follows. The elementary definitions, and the symmetric operations, 
is treated in {\bf section $2$}, in {\bf section $3$} we give the two relations between 
avoidance of patterns in $S_k$ and avoidance of signed patterns in $E_k^r$, in {\bf section $4$} 
we represent two sets of signed patterns, and represent a bijection which gives 
a combinatorial geometric explanation for one of these results. 
In {\bf sections $5$, $6$} we prove the first and the second part of Main Theorem, 
respectively. Finally, in the {\bf last section} we prove a combinatorial identity 
as a corollary of Main Theorem. \\
\\
{\bf Main Theorem:} \\
$(i)$ For every $2$-letter signed pattern $[\tau]_a$, the
number of $[\tau]_a$-avoiding signed permutations in $E_n^r$ is given by 
the expression: $\sum\limits_{j=0}^n j!(r-1)^j{n\choose j}^2$.
\\
$(ii)$ A double restriction by $2$-letter signed patterns gives  
one $1$-Wilf class, four $2$-Wilf clases, six $r$-Wilf clases for $r\geq 3$.
%============================================================================
\section*{\centering {\sc 2. Symmetries on signed permutations}}

As on the symmetric group $S_n$ there are two natural symmetric operations, 
the reversal and the complement (see ~\cite{SS}), also on $E_n^r$ we define:
\begin{enumerate}
\item[(i)] the {\em reversal} $er:E_n^r\rightarrow E_n^r$ defined by 
   $$er:(\alpha_1^{(u_1)},\dots,\alpha_n^{(u_n)})\mapsto  (\alpha_n^{(u_n)},\dots,\alpha_1^{(u_1)});$$

\item[(ii)] the {\em complement} $ec:E_n^r\rightarrow E_n^r$ defined by
   $$ec:(\alpha_1^{(u_1)},\dots,\alpha_n^{(u_n)})\mapsto ((n+1-\alpha_1)^{(u_1)},\dots,(n+1-\alpha_n)^{(u_n)});$$

\item[(iii)] and besides that, there is the {\em sign-complement} $es:E_n^r\rightarrow E_n^r$ 
defined by 
   $$es:(\alpha_1^{(u_1)},\dots,\alpha_n^{(u_n)})\mapsto (\alpha_1^{(r+1-u_1)},\dots,\alpha_n^{(r+1-u_n)}).$$
\end{enumerate}

\begin{example}
Let $\Phi=(1^{(1)},3^{(2)},2^{(1)})\in E_3^2$, then 
$er(\Phi)=(2^{(1)},3^{(2)},1^{(1)})$, $ec(\Phi)=(2^{(1)},1^{(2)},3^{(1)})$, and 
$es(\Phi)=(1^{(2)},2^{(1)},3^{(2)})$.
\end{example}

\begin{proposition}
The group $<er,ec,es>$ is isomorphic to $D_8$.
\end{proposition}

More generally, we extend these symmetric operations to subsets of $E_n^r$:
$g(T)=\{g(\Phi)|\Phi\in T\}$, where $g=er$, $ec$, or $es$. 

\begin{theorem}
\label{rcq}
Let $T\subset E_k^r$. For all $n\geq 0$, 
	$$|E_n^r(T)|=|E_n^r(er(T))|=|E_n^r(ec(T))|=|E_n^r(es(T))|.$$
\end{theorem}

Now we define the fourth symmetric operation on $E_n^r$. Let us define
	$$h_{\delta,n}:E_n^r\rightarrow E_n^r,$$ 
where $\delta\in S_r$ by $h_{\delta,n}([\alpha]_a)=[\alpha]_b$
such that $b_i=\delta_{a_i}$ for all $i=1,2,\dots,n$. More generally,  
  $h_{\delta,n}(T)=\{h_{\delta,n}([\alpha]_a) | [\alpha]_a\in T\}$
for $T\subset E_n^r$. 

\begin{theorem}
\label{prem}
Let $T\subset E_k^r$, $\delta\in S_r$. Then $|E_n^r(T)|=|E_n^r(h_{\delta,k}(T))|$.
\end{theorem}
\begin{proof}

Let $[\alpha]_a\in E_n^r(T)$, so $[\alpha]_a$ is $T$-avoiding if and 
only if $h_{\delta,n}([\alpha]_a)$ is $h_{\delta,k}(T)$-avoiding. On 
the other hand $h_{\delta,n}$ is an invertible function. Hence 
the theorem holds.
\qed\end{proof}

\begin{corollary}
\label{lxx}
Let $T\subseteq E_k^r$, and let $\delta\in S_r$ such that $a_{b_j}=j$ for $j=1,2,\dots,d$. 
For all $n\geq 0$, $|E_n^r(T)|=|E_n^r(h_{\delta,k}(T))|$.
\end{corollary}

\begin{example}
As an example, for $r\geq 3$,
 $$|E_n^r((1^{(1)},2^{(2)}),(1^{(2)},2^{(3)}))|=|E_n^r((1^{(2)},2^{(1)}),(1^{(1)},2^{(3)}))|,$$ 
by the symmetric operation $h_{(2,1,3,4,\dots,r),n}$.
\end{example}
%=========================================================================
\section*{\centering {\sc 3. Avoidance patterns and signed patterns}}

We say a signed permutation $[\tau]_a\in E_k^r$ is {\em homogeneous} 
if $a_i=u$ for all $i=1,2,\dots,k$ where $1\leq u\leq r$; in this case 
we denote $[\tau]_a$ by $[\tau]_{(u)}$. More generally, we denote 
$T_{(u)}=\{[\tau]_{(u)}| \tau\in T\}$.

\begin{theorem}
\label{lsg1}
Let $1\leq u\leq r$, $T\subset S_k$. For all $n\geq 0$
    $$|E_n^r(T_{(u)})|=\sum_{j=0}^n j!(r-1)^j|S_{n-j}(T)| {n\choose j}^2.$$
\end{theorem}
\begin{proof}

Immediately by definitions 
$$|E_n^r(T_{(u)})|=\sum\limits_{j=0}^n {n\choose j}^2 |E_{\{1,2,\dots,j\}}^{\{u\}}(T_{(u)})| 
		|E_{\{j+1,\dots,n\}}^{\{1,\dots,u-1,u+1,\dots,r\}}|,$$
where $E_{T_1}^{T_2}$ is the set of all signed 
permutations with set symbols $T_1$ and set signs $T_2$. 
So clearly $|E_{\{j+1,\dots,n\}}^{\{1,\dots,u-1,u+1,\dots,r\}}|=(n-j)!\cdot (r-1)^{n-j}$, 
also $|E_{\{1,2,\dots,j\}}^{\{u\}}(T_{(u)})|=|S_j(T)|$ by removing the sign $u$. Hence the 
theorem holds.
\qed\end{proof} 

\begin{example} ({\em see} \cite[Eq. $46$]{Rs})
 For $a=1,2$, by Theorem \ref{lsg1}, 
	$$|E_n^2((12)_{(a)},(21)_{(a)} )|=(n+1)!,$$
	$$|E_n^2( (12)_{(a)} )|=|E_n^2( (21)_{(a)} )|=\sum\limits_{j=0}^n j!{n\choose j}^2.$$ 
\end{example}

\begin{theorem}
\label{superset}
Let $r\geq 1$, $\tau\in S_k$. For all $n\geq 0$, $|E_n^r(F_{\tau})|=r^n |S_n(\tau)|$, 
where $F_{\tau}=\{(\tau_1^{(v_1)},\dots,\tau_k^{(v_k)}) | 1\leq v_i\leq r \}$.
\end{theorem}
\begin{proof}

Let us define a function $f:[r]^n\times S_n(\tau)\mapsto E_n^r(F_{\tau})$ by 
 $$f((u_1,\dots,u_n;\alpha_1,\dots,\alpha_n))=(\alpha_1^{(u_1)},\dots,\alpha_n^{(u_n)}).$$
So $(u_1,\dots,u_n;\alpha_1,\dots,\alpha_n)\in [r]^n\times S_n(\tau)$  
if and only if $(\alpha_1,\dots,\alpha_n)$ avoids $\tau$, which is equivalent to 
$(\alpha_1^{(u_1)},\dots,\alpha_n^{(u_n)})$ avoids $F_{\tau}$ for all 
$u_i=1,2,\dots,r$. Hence $f$ is a bijection, which means that the theorem holds.
\qed\end{proof}

\begin{example}
Let $T=\{ (1^{(a)},2^{(b)})| a,b=1,2,\dots,r \}$; by 
Theorem \ref{superset} we obtain $|E_n^r(T)|=r^n$ for all $n\geq 0$.
\end{example}
%=====================================================================================
\section*{\centering {\sc 4. Restricted sets}}

In this section, we calculate cardinalities of $E_n^r(T)$ 
for two special subsets $T\subset E_2^r$. The first special subset is
defined by $T_{b;a_1,a_2,\dots,a_l}=\{ (1^b,2^{(a_j)}) | j=1,2,\dots,l\}$.

\begin{theorem}
\label{l12l}
Let $1\leq l \leq r$, and $1\leq b\leq a_1<a_2<\dots<a_l\leq r$. Then
$$\sum\limits_{n\geq 0} \frac{|E_n^r(T_{b;a_1,a_2,\dots,a_l} )|}{n!} x^n=
	\left( \frac{1-(r-l)x}{( 1-(r-1)x)^l} \right) ^{\frac1{l-1}};$$
when $l=1$ we take the limit of the right hand side which 
equals $\frac{e^{\frac{1}{1-(r-1)x}}} {1-(r-1)x}$.
\end{theorem}
\begin{proof}
By Corollary \ref{lxx} $|E_n^r(T_{b;a_1,\dots,a_l})|=|E_n^r(1;a,a+1,\dots,a+l-1)|$. 

Let $\Phi\in E_n^r(T_{1;a+1,\dots,a+l-1})$, $p_r(n)=|E_n^r(T_{1;a+1,\dots,a+l-1})|$,  
and let us consider the possible values of $\Phi_1$: 
\begin{enumerate}
\item	Let $\Phi_1=i^{(c)}$, $c\neq 1$, and $1\leq i\leq n$; so 
	$\Phi\in E_n^r(T_{1;a+1,\dots,a+l-1})$ 
	if and only if $(\Phi_2,\dots,\Phi_n)$ is $T_{1;a+1,\dots,a+l-1}$
	-avoiding, hence in this case there are $(r-1)np_r(n-1)$ 
	signed permutations. 

\item	Let $\Phi_1=i^{(1)}$; since $\Phi$ is $T_{1;a+1,\dots,a+l-1}$-avoiding, 
	the symbols $i+1,\dots,n$ appeared with sign $d\geq a+1$ or $d\leq a-1$. 
	Also the symbols $1,\dots,i-1$ are $T_{1;a+1,\dots,a+l-1}$-avoiding, and
          can be replaced anywhere at positions 
	$2,\dots,n$, hence there are 
	$\sum\limits_{i=1}^n {{n-1}\choose{i-1}} |E_{\{i+1,\dots,n\}}^{\{1,\dots,a-1,a+l,\dots,r\}}|\cdot |E_{i-1}^r(T)|$ 
	signed permutations, which means there are  
	$\sum\limits_{i=1}^n {{n-1}\choose{i-1}} (n-i)!(r-l)^{n-i}p_r(i-1)$ 
	signed permutations.
\end{enumerate}		
So by the above two cases we obtain a recurrence relation satisfied by
$p_n$ 
  $$p_n=(r-1)np_{n-1}+\sum_{i=1}^n {{n-1}\choose{i-1}} (n-i)!(r-l)^{n-i} p_{i-1},$$
for $n\geq 1$, and $p_0=1$. Let $q_n=p_r(n)/n!$. By multiplying the recurrence 
by $x^{n-1}/(n-1)!$, and summing up over all $n\geq 1$, we obtain 
  $$\frac{d}{dx} q(x)=(r-1)\frac{d}{dx} (x q(x))+\frac{q(x)}{1-(r-l)x},$$
where $q(x)$ is the generating function of $q_n$. Besides $q(0)=1$,
hence the theorem holds.
\qed\end{proof}

\begin{corollary}
For all $n\geq 0$, $|E_n^r(T_{1;1,2,\dots,r})|=\prod\limits_{j=0}^n (j(r-1)+1).$
\end{corollary}
\begin{proof}
Immediately by the proof of Theorem \ref{l12l}, for $n\geq 2$
   $$p_r(n)=|E_n^r(T_{1;1,2,\dots,r})|=((r-1)n+1)p_r(n-1).$$
Besides, $p_r(1)=r$, and $p_r(0)=1$, hence the corollary holds.
\qed\end{proof}

\begin{example} ({\em see} \cite[Eq. $47$]{Rs})
By Theorem \ref{l12l}, $|E_n^2((1^{(1)},2^{(1)}),(1^{(1)},2^{(2)}))|=(n+1)!$.
\end{example}

Now we represent the second special subset. 
Consider a subset $T\subset E_k^r$; we say that $T$ is {\em good} if it is the union of
disjoint homogeneous subsets; that is, $T=\cup_{j=1}^p (T_j)_{(u_j)}$.
As an example, $T=\{123_{(1)},132_{(1)},213_{(2)}\}$ is a good set. 

\begin{theorem}
\label{disgood}
Let $T=\cup_{j=1}^p (T_j)_{(u_j)}$ be a good set. For $n\geq 0$
    $$|E_n^r(T)|=\sum_{j_1=0}^n\sum_{j_2=0}^{n-j_1}\dots\sum_{j_p=0}^{n-j_1-\dots -j_{p-1}} (r-p)^{n-j_1-\dots -j_p} 
		\frac{{n\choose {j_1,j_2,\dots,j_p}}^2}{(n-j_1-\dots -j_p)!} \prod_{i=1}^p |S_{j_i}(T_i)|.$$
\end{theorem}
\begin{proof}

The theorem holds for $p=1$ by Theorem \ref{lsg1}. Now let $p>1$, so by definitions
    $|E_n^r(T)|=\sum_{j_1=0}^n |E_{n-j}^{1,\dots,u_1-1,u_1+1,\dots,r\}}(T\backslash (T_1)_{(u_1)})| |S_{j}(T_1)| {n\choose {j_1}}^2,$
therefore, 
    $|E_n^r(T)|=\sum\limits_{j_1=0}^n |E_{n-j}^{r-1}(T\backslash (T_1)_{(u_1)})| |S_j(T_1)| {n\choose {j_1}}^2.$
Hence by induction the theorem holds. 
\qed\end{proof}

Let $T_{d,l;a_1,\dots,a_l}$ be a subset of $E_2^k$ defined by 
	$$T_{d,l;a_1,\dots,a_l}=
         \cup_{i=1}^d \{(1,2)_{(a_i)}\} \bigcup \cup_{i=d+1}^l \{(2,1)_{(a_i)}\},$$
hence by Theorem \ref{disgood} we obtain the following corollary:

\begin{corollary}
\label{ll12l}
Let $1\leq a_1,\dots,a_l\leq k$ be $l$ different numbers. For $n\geq 0$, 
  $$|E_n^r(T_{d,l;a_1,\dots,a_l})|=
   \sum_{i_1+\dots+i_l\leq n} \frac{{n\choose {i_1,\dots,i_l}}^2}{(n-i_1-\dots -i_l)!} (r-l)^{n-i_1-\dots-i_l}.$$
\end{corollary}

Now we built a bijection, which gives for the set 
$E_n^r(T_{d,a;a_1,\dots,a_l})$ a combinatorial geometric explanation. 
Consider $l$ lines $L_1$,\dots,$L_l$ such that $L_i$ contains 
all the points of the form $j^{(i)}$ for all $j=1,2,\dots,n$. We say 
$L_i$ is {\em good} if the points $1^{(i)}$ to $n^{(i)}$ are decreasing, 
and the line $L_i$ is {\em bad} if the points $1^{(i)},\dots,n^{(i)}$ are 
increasing, otherwise we say the line $L_i$ is {\em free}.\\

Now we consider the following collection which represents the set 
$T_{d,l;a_1,\dots,d_l}$. Let $L_{a_1},\dots,L_{a_d}$ be good lines,  
$L_{a_{d+1}},\dots,L_{a_l}$ be bad lines, and $L_i$ be a free line for all 
$1\leq i\leq k$ such that $i\not\in\{a_1,\dots,a_l\}$. For example,  
the representation of $T_{1,2;3,2}$ where $k=4$, is given by the following 
diagram.
	\begin{center}
	\unitlength=.2mm
	\begin{picture}(500,100)(0,-20)
	\put(111,80){\line(1,-1){80}}
	\put(98,80){\makebox(0,0){$n^{^{(3)}}$}}
	\put(154,20){\makebox(0,0){$3^{^{(3)}}$}}
	\put(164,10){\makebox(0,0){$2^{^{(3)}}$}}
	\put(174,0){\makebox(0,0){$1^{^{(3)}}$}}
	\put(201,-12){\makebox(0,0){$^{L_3}$}}

	\put(183,80){\line(1,-1){80}}
	\put(168,80){\makebox(0,0){$1^{^{(2)}}$}}
	\put(178,70){\makebox(0,0){$2^{^{(2)}}$}}
	\put(188,60){\makebox(0,0){$3^{^{(2)}}$}}
	\put(248,0){\makebox(0,0){$n^{^{(2)}}$}}
	\put(271,-12){\makebox(0,0){$^{L_2}$}}

	\put(242,80){\line(1,-1){80}}
	\put(331,-12){\makebox(0,0){$^{L_1}$}}

	\put(302,80){\line(1,-1){80}}
	\put(391,-12){\makebox(0,0){$^{L_4}$}}
	\end{picture}
	\end{center}
\begin{center} 	{{\bf Figure 1}:\ Representation of $T_{1,2;3,2}$} \end{center}
Here the lines $L_1$ and $L_4$ are free lines.\\

Now let us define a {\em path} between the points on the lines of the 
representation of $T_{d,l;a_1,\dots,a_l}$. A path is a collection of steps, 
starting anywhere, such that every step is one of the following 
steps: 

\begin{enumerate}
\item[(i)] a decreasing step from a point to another point on a bad, or a good line,

\item[(ii)] a free step on the free line, or between the lines
(from a point to another point).
\end{enumerate}

Hence by definitions we immediately have the following proposition.

\begin{proposition}
\label{pp1}
Every path of $n$ steps is a $T_{d,l;a_1,\dots,a_l}$-avoiding 
signed permutation in $E_n^r$. 
\end{proposition}

Using the above proposition we find the cardinality of the set 
$E_n^r(T_{d,l;a_1,\dots,a_l})$ by the following theorem.

\begin{theorem}
\label{lines12l}
Let $a_1,\dots,a_l$ be $l$ different numbers such that $1\leq a_i\leq r$ 
for all $i=1,2,\dots,l$. For $n\geq 0$,
    $$|E_n^r(T_{d,l;a_1,\dots,a_l})|= 
	\sum_{i_1+\dots+i_l\leq n} \frac{{n\choose {i_1,\dots,i_l}}^2}{(n-i_1-\dots-i_l)!}(r-l)^{n-i_1-\dots-i_l}.$$
\end{theorem}
\begin{proof}
To choose a path of $n$ steps with $l$ points in bad or good lines we have to:
\begin{enumerate}
\item[(i)]	choose $i_1,\dots i_l$ places in the path. There are 
	${n\choose {i_1,\dots,i_l}}$ possibilities.

\item[(ii)]	choose ${i_1,\dots,i_l}$ points from bad or good  
	lines. There are 
${n\choose {i_1,\dots,i_l}}$ possibilities.

\item[(iii)]	choose $n-d$ points on free lines, where 
	$d=i_1+\dots+i_l$. There are $(n-d)!(k-l)^{n-d}$ possibilities.
\end{enumerate}
Hence, by Proposition \ref{pp1} the theorem holds
\qed\end{proof}

By Theorem \ref{lines12l} we obtain a generalization of certain results in ~\cite{Rs}, 
particularly we get the following corollary.

\begin{corollary}
\label{c12l}
Let $0\leq d\leq r$; for $n\geq 0$,
   $$|E_n^r(T_{d,r;1,2,3,\dots,r})|=\sum_{i_1=0}^n\sum_{i_2=0}^{n-i_1}\dots\sum_{i_{r-1}=0}^{n-i_1-\dots,i_{r-2}} {n\choose {i_1,\dots,i_{r-1},n-i_1-\dots-i_{r-1}}}^2.$$
\end{corollary}
 
\begin{example} ({\em see} \cite[Eq. $49$]{Rs}) 
By Corollary \ref{c12l} we obtain for $\beta,\gamma\in S_2$
	$$|E_n^2(\beta_{(1)},\gamma_{(2)})|=\sum_{i=0}^n {n\choose i}^2={{2n}\choose n}.$$ 
\end{example}
%===========================================================================
\section*{\centering{\sc 5. Single restriction by a $2$-letter signed pattern}} 

The length $2$ signed permutations give rise to some enumeratively
interesting classes of signed permutations, which we examine in this
section.
In the symmetric group $S_n$, patterns of length $2$ are uninterestingly
restrictive, and length $3$ is the first interesting case. Also in
$E_n^r$, restriction by patterns of length $1$ is trivial, and given by the 
following formula $|E_n^r(1^a)|=n!\cdot(r-1)^n$, where $1\leq a\leq r$.
Let us denote $d_r(n)=\sum\limits_{j=0}^n j!(r-1)^j{n\choose j}^2$, and let $d_r(x)$ 
be the generating function of the sequence $d_r(n)/n!$. Hence it easy to see that 
$d_r(x)=\frac{e^\frac{x}{1-(r-1)x}}{1-(r-1)x}.$\\

Now we prove the first case of Main Theorem, that is, that there exists exactly one $r$-Wilf class 
of a single restriction by a $2$-letter signed pattern, for all $r\geq 1$.

\begin{theorem}
\label{oneab}
Let $r\geq 1$, and $1\leq a,b,c,d\leq r$. For $n\geq 0$
	  $$|E_n^r((1^{(a)},2^{(b)}))|=|E_n^r((2^{(c)},1^{(d)}))|=d_r(n).$$
\end{theorem}
\begin{proof}
By section $2$ (symmetric operations) we have to prove the following 
two cases:
\begin{enumerate}
\item	Let $1\leq a\leq r$; for $n\geq 0$, $|E_n^r((1^{(a)},2^{(a)}))|=|E_n^r((2^{(a)},1^{(a)}))|=d_r(n)$; 

\item	Let $b\leq a$; for $n\geq 0$, $|E_n^r((1^{(a)},2^{(b)}))|=|E_n^r((1^{(a)},2^{(a)}))|$. 
\end{enumerate}

The first, and the second cases are obtained immediately by Theorem \ref{lsg1}, and 
Theorem \ref{l12l}, respectively. 
\qed\end{proof}
%==========================================================================
\section*{\centering{\sc 6. Double restrictions by $2$-letter signed patterns}}

In this section, we find the number of $r$-Wilf clases, $r\geq 1$, of double 
restrictions by $2$-letter signed patterns.
In $E_2^r$ there are $k^2(k^2-1)$ possibilities to choose two elements of the following form:
$(1^{(a)},2^{(b)}), (1^{(c)},2^{(d)}),$
and there are $k^4$ possibilities to choose two lements of the following form:
$(1^{(a)},2^{(b)}), (2^{(c)},1^{(d)}),$
where $1\leq a,b,c,d\leq r$. On the other hand, by symmetric operations 
(section 2), the question of determining the 
$E_n^r([\tau]_a, [\tau ']_{a'})$ for $k^2(2k^2-1)$ choices for $2$-letter signed patterns 
$[\tau]_a$, $[\tau']_{a'}$ reduces to determining the 
$E_n^r([\tau]_a, [\tau ']_{a'})$ where $[\tau]_a$, $[\tau']_{a'}$ are from Table \ref{tab1}.

\begin{table}
\begin{center}
    \begin{tabular}{|c|c|c|c|c|} \hline
Case  &	$[\tau]_a$ 			& $[\tau']_{a'}$	& $|E_n^5([\tau]_a,[\tau']_{a'})|$ for $n=0,1,2,3,4,5$ & Reference	\\ \hline\hline
1&\	$(1^{(1)},2^{(1)})$\ 		&\ $(2^{(1)},1^{(1)})$\	& $1,5,48,672,12288,276480$	& Theorem \ref{case1}			\\ \hline %aa
2&\	$(1^{(1)},2^{(1)})$\		&\ $(1^{(1)},2^{(2)})$\	& $1,5,48,672,12288,276480$	& Theorem \ref{case1}			\\ \hline %aa
3&\	$(1^{(1)},2^{(2)})$\		&\ $(2^{(1)},1^{(2)})$\	& $1,5,48,672,12288,276480$	& Theorem \ref{case1}			\\ \hline %aa
4&\	$(1^{(1)},2^{(2)})$\		&\ $(2^{(2)},1^{(1)})$\	& $1,5,48,672,12288,276480$	& Theorem \ref{case1}			\\ \hline %aa
5&\	$(1^{(1)},2^{(2)})$\		&\ $(1^{(1)},2^{(3)})$\	& $1,5,48,672,12288,276480$	& Theorem \ref{case1}			\\ \hline\hline %aa

6&\	$(1^{(1)},2^{(1)})$\ 		&\ $(1^{(2)},2^{(2)})$\	& $1,5,48,668,12046,265062$	& Theorem \ref{case2}			\\ \hline %aa
7&\	$(1^{(1)},2^{(1)})$\ 		&\ $(2^{(2)},1^{(2)})$\	& $1,5,48,668,12046,265062$	& Theorem \ref{case2}			\\ \hline %aa	
8&\	$(1^{(1)},2^{(1)})$\		&\ $(1^{(2)},2^{(3)})$\	& $1,5,48,668,12046,265062$	& Theorem \ref{case2}			\\ \hline %aa
9&\	$(1^{(1)},2^{(1)})$\		&\ $(2^{(2)},1^{(3)})$\	& $1,5,48,668,12046,265062$	& Theorem \ref{case2}			\\ \hline %aa
10&\	$(1^{(1)},2^{(2)})$\		&\ $(1^{(3)},2^{(4)})$\	& $1,5,48,668,12046,265062$	& Theorem \ref{case2}			\\ \hline %aa
11&\	$(1^{(1)},2^{(2)})$\		&\ $(2^{(3)},1^{(4)})$\	& $1,5,48,668,12046,265062$	& Theorem \ref{case2}			\\ \hline\hline %aa

12&\	$(1^{(1)},2^{(2)})$\		&\ $(2^{(1)},1^{(3)})$\	& $1,5,48,670,12168,270856$	& Theorem \ref{case3}			\\ \hline %aa
13&\	$(1^{(1)},2^{(2)})$\		&\ $(2^{(2)},1^{(3)})$\	& $1,5,48,670,12168,270856$	& Theorem \ref{case3}			\\ \hline %aa
14&\	$(1^{(1)},2^{(2)})$\		&\ $(2^{(3)},1^{(1)})$\	& $1,5,48,670,12168,270856$	& Theorem \ref{case3}			\\ \hline\hline %aa

15&\	$(1^{(1)},2^{(1)})$\		&\ $(2^{(1)},1^{(2)})$\	& $1,5,48,671,12288,273665$	& Theorem \ref{case4}			\\ \hline\hline %aa

16&\	$(1^{(1)},2^{(2)})$\		&\ $(1^{(2)},2^{(3)})$\	& $1,5,48,669,12106,267867$	& 					\\ \hline\hline %aa

17&\	$(1^{(1)},2^{(2)})$\		&\ $(1^{(2)},2^{(1)})$\	& $1,5,48,670,12166,270672$	& 					\\ \hline
    \end{tabular}
\end{center}
   \caption {Pairs of $2$-letter signed patterns}
   \label{tab1}
\end{table}
\begin{theorem}
\label{case1}
For $n\geq 0$, $|E_n^r(T)|=n!(n+r-1)(r-1)^{n-1}$ where 
\begin{enumerate}
\item[(i)]	$T=\{ (1^{(1)},2^{(1)}), (2^{(1)},1^{(1)}) \}$ for $r\geq 1$;

\item[(ii)] 	$T=\{ (1^{(1)},2^{(1)}), (1^{(1)},2^{(2)}) \}$ for $r\geq 2$;
 
\item[(iii)] 	$T=\{ (1^{(1)},2^{(2)}), (2^{(1)},1^{(2)}) \}$ for $r\geq 2$;

\item[(iv)] 	$T=\{ (1^{(1)},2^{(2)}), (2^{(2)},1^{(1)}) \}$ for $r\geq 2$;

\item[(v)] 	$T=\{ (1^{(1)},2^{(2)}), (1^{(1)},2^{(3)}) \}$ for $r\geq 3$.
\end{enumerate}
\end{theorem}
\begin{proof}
By Theorem \ref{lsg1} it is easy to obtain $(i)$, and Theorem \ref{l12l} immediately 
yields $(ii)$, and $(v)$. Now let us prove $(iii)$ and $(iv)$.\\

{\bf Case $(iii)$:} Let $p_n=|E_n^r(T)|$, $\Phi\in E_n^r(T)$, 
and let us consider the possible values of $\Phi_1$:
\begin{enumerate}
\item	Let $\Phi_1=i^{(c)}$, $c\neq 1$; $\Phi\in E_n^r(T)$
	if and only if $(\Phi_2,\dots,\Phi_n)\in E_{\{1,\dots,i-1,i+1,\dots,n\}}^r(T)$. 
	Hence in this case there are $(r-1)np_{n-1}$ signed permutations.

\item	Let $\Phi_1=i^{(1)}$; since $\Phi$ is $T$-avoiding, 
	the symbols $1,\dots,i-1,i+1,\dots,n$ are not signed by $2$, and can be replaced 
	anywhere at positions $2,\dots,n$. Hence, in this case there are $(n-1)!(r-1)^{n-1}$  
	signed permutations.
\end{enumerate}
Therefore by the above three cases $p_n$ staisfies the following relation:
	  $$p_n=n(r-1)p_{n-1}+n!(r-1)^{n-1}.$$
Besides $p_0=1$, and $p_1=r$, hence $(iv)$ holds.\\

{\bf Case $(iv)$:} Let $p_n=|E_n^r(T)|$, $\Phi\in E_n^r(T)$ such that $\Phi_j=n^{(c)}$, 
and let us consider the possible values of $j,c$:
\begin{enumerate}
\item	Let $c\neq 2$; $\Phi\in E_n^r(T)$
	if and only if $(\Phi_1,\dots,\Phi_{j-1},\Phi_{j+1},\dots,\Phi_n)\in E_{n-1}^r(T)$. 
	Hence in this case there are $(r-1)np_{n-1}$ signed permutations.

\item	Let $c=2$; $\Phi\in E_n^r(T)$
	if and only if $(\Phi_1,\dots,\Phi_{j-1},\Phi_{j+1},\dots,\Phi_n)$ is a signed 
	permutation with symbols $1,2,\dots,n-1$ and signs $2,\dots,r$. 
	Hence, in this case there are $(n-1)!(r-1)^{n-1}$ signed permutations.
\end{enumerate}
Therefore by the above three cases $p_n$ staisfies the following relation:
	  $$p_n=n(r-1)p_{n-1}+n!(r-1)^{n-1}.$$
Besides $p_0=1$, and $p_1=r$, hence $(v)$ holds.
\qed\end{proof}

\begin{example} ({\em see} \cite[Eq. $46$, $47$]{Rs})
 As an example we get 
$$\begin{array}{ll}
|E_n^2((1^{(1)},2^{(1)}),(2^{(1)},1^{(1)}))|=|E_n^2((1^{(1)},2^{(1)}),(1^{(1)},2^{(2)}))|=& \\
|E_n^2((1^{(1)},2^{(2)}),(2^{(1)},1^{(2)}))|=|E_n^2((1^{(1)},2^{(2)}),(2^{(2)},1^{(1)}))|=& (n+1)!
\end{array}$$ for $n\geq 0$, 
which was proved in \cite{Rs}.
\end{example}
%====================================================
%====================================================
\begin{theorem}
\label{case2}
Let $2\leq a\leq b$, and $r\geq b$; for all $n\geq 1$
  $$|E_n^r(T)|=\sum\limits_{i+j\leq n} {n\choose {i,j,n-i-j}}^2 (n-i-j)!(r-2)^{n-i-j},$$
where 
\begin{enumerate}
\item[(i)] 	$T=\{ (1^{(1)},2^{(1)}),(1^{(a)},2^{(b)})\}$; 

\item[(ii)]	$T=\{ (1^{(1)},2^{(1)}),(2^{(a)},1^{(b)})\}$;

\item[(iii)]	$T=\{ (1^{(1)},2^{(2)}),(1^{(3)},2^{(4)})\}$;

\item[(iv)]	$T=\{ (1^{(1)},2^{(2)}),(2^{(3)},1^{(4)})\}$.
\end{enumerate}
\end{theorem}
\begin{proof}
{\bf Cases $(i)$, $(ii)$:} Immediately by the proof of Theorem \ref{disgood}, and by 
part $(i)$ of Main Theorem, we claim these cases.\\

{\bf Cases $(iii)$, $(iv)$:\ } Let $T_1=\{ (1^{(1)},2^{(2)}), (1^{(3)},2^{(4)}) \}$, $T_2=\{ (1^{(1)},2^{(2)}), (1^{(3)},2^{(4)}) \}$,
and let $\Phi\in E_n^r(T_1)$. Also let us define $I_{\Phi}$ to be the set of all $j$ such that 
$\Phi_j$ is signed by either $3$ or $4$. Now we define a function $f:E_n^r(T_1)\rightarrow E_n^r(T_2)$ 
by reversing all the $\Phi_j$ where $j\in I_{\Phi}$. Hence by definitions, $f$ is a bijection, which means that 
$|E_n^r((1^{(1)},2^{(2)}),(2^{(3)},1^{(4)}))|=|E_n^r((1^{(1)},2^{(2)}),(2^{(4)},1^{(3)}))|$. \\

On the other hand, immediately by Main Theorem part $(i)$ we get 
$$|E_n^r((1^{(1)},2^{(2)}),(2^{(3)},1^{(4)}))|=|E_n^r((1^{(1)},2^{(1)}),(2^{(4)},1^{(4)}))|,$$ 
which means by case $(i)$ that the theorem holds.
\qed\end{proof}

\begin{example} ({\em see} \cite[Eq. $47$]{Rs})
As an example, by Theorem \ref{case2} for $n\geq 0$
	$$|E_n^2((1^{(1)},2^{(1)}),(2^{(2)},1^{(2)}))|={{2n}\choose n}.$$  
\end{example}
%=================================================================
%=================================================================
\begin{theorem}
\label{case3}
For $r\geq 3$, $\sum\limits_{n\geq 0} \frac{E_n^r(T)}{n!}x^n=\frac{\int d_{r-1}^2(x) dx}{1-(r-1)x}$ 
where
\begin{enumerate}
\item[(i)]	$T=\{ (1^{(1)},2^{(2)}), (2^{(1)},1^{(3)}) \}$; 

\item[(ii)]	$T=\{ (1^{(1)},2^{(2)}), (2^{(2)},1^{(3)}) \}$;

\item[(iii)]	$T=\{ (1^{(1)},2^{(2)}), (2^{(3)},1^{(1)}) \}$.
\end{enumerate}
\end{theorem}
\begin{proof}

{\bf Case $(i)$:\ } Let $T=\{(1^{(1)},2^{(2)}),(2^{(1)},1^{(3)}))\}$, 
$p_n=E_n^r(T)$, $\Phi\in E_n^r(T)$, and let us consider the possible values of $\Phi_1$:
\begin{enumerate}
\item 	Let $\Phi_1=i^{(c)}$, $c\neq 1$; so $\Phi\in E_n^r(T)$
	if and only if $(\Phi_2,\dots,\Phi_n)\in E_{\{1,\dots,i-1,i+1,\dots,n\}}^r(T)$. 
	Hence in this case there are $(r-1)np_{n-1}$ signed permutations.

\item	Let $\Phi_1=i^{(1)}$; since $\Phi$ is $T$-avoiding,   
	the symbols $i+1,\dots,n$ are not signed by $2$, and the symbols 
	$1,\dots,i-1$ are not signed by $3$. Hence there are 
	${{n-1}\choose {i-1}} |E_{n-i}^{r-1}((2^{(1)},1^{(3)}))| |E_{i-1}^{r-1}((1^{(1)},2^{(2)}))|$ 
	signed permutations, which means by part $(i)$ of Main Theorem that there are  
	${{n-1}\choose {i-1}}d_{r-1}(n-i)d_{r-1}(i-1)$ signed permutations.
\end{enumerate}
Therefore $p_n$ staisfies the following relation:
   $$p_n=n(r-1)p_{n-1}+\sum_{i=1}^n {{n-1}\choose {i-1}} d_{r-1}(n-i) d_{r-1}(i-1).$$
Besides $p_0=1$, and $p_1=r$, hence case $(i)$ holds.\\

{\bf Case $(ii)$:} Let $T=\{(1^{(1)},2^{(2)}),(2^{(2)},1^{(3)}))\}$, $p_n=E_n^r(T)$, and 
$\Phi\in E_n^r(T)$ such that $\Phi_j=n^{(c)}$. Let us consider  
the possible values of $j,c$:
\begin{enumerate}
\item	Let $c\neq 2$; $\Phi\in E_n^r(T)$
	if and only if $(\Phi_1,\dots,\Phi_{j-1},\Phi_{j+1},\dots,\Phi_n)\in E_{n-1}^r(T)$. 
	Hence in this case there are $(r-1)np_{n-1}$ signed permutations.

\item	Let $c=2$; since $\Phi$ is $T$-avoiding,  
	all the symbols in $(\Phi_1,\dots,\Phi_{j-1})$ are not signed by $1$, and the symbols in 
	$(\Phi_{j+1},\dots,\Phi_n)$ are not signed by $3$. 
	Hence there are ${{n-1}\choose {j-1}} |E_{n-j}^{r-1}((2^{(2)},1^{(3)}))| |E_{j-1}^{r-1}((1^{(1)},2^{(2)}))|$ 
	signed permutations, which means by part $(i)$ of Main Theorem that there are 
	${{n-1}\choose {j-1}} d_{r-1}(n-j)d_{r-1}(j-1)$ signed permutations.
\end{enumerate}
Therefore $p_n$ staisfies the following relation:
   $$p_n=n(r-1)p_{n-1}+\sum_{j=1}^n {{n-1}\choose {j-1}} d_{r-1}(n-j) d_{r-1}(j-1).$$
Besides $p_0=1$, and $p_1=r$, hence case $(ii)$ holds.\\

{\bf Case $(iii)$:} Similarly to the case $(ii)$ for 
$\Phi\in E_n^r(T)$ such that $\Phi_j=1^{(c)}$, we consider the possible 
values of $j,c$, and get the same result.
\qed\end{proof}
%=================================================================
%=================================================================
\begin{theorem}
\label{case4}
For $r\geq 2$, $\sum\limits_{n\geq 0} \frac{E_n^r( (1^{(1)},2^{(1)}), (2^{(1)},1^{(2)}) )}{n!}x^n=
		\frac{\int \frac{d_{r-1}(x)}{1-(r-1)x}dx}{1-(r-1)x}.$
\end{theorem}
\begin{proof}

Let $T=\{(1^{(1)},2^{(1)}),(2^{(1)},1^{(2)}))\}$, $p_n=E_n^r(T)$, $\Phi\in E_n^r(T)$, 
and let us consider the possible values of $\Phi_1$:
\begin{enumerate}
\item	If $\Phi_1=i^{(c)}$ where $c\neq 1$, then $\Phi\in E_n^r(T)$
	if and only if $(\Phi_2,\dots,\Phi_n)\in E_{\{1,\dots,i-1,i+1,\dots,n\}}^r(T)$. 
	Hence in this case there are $(r-1)np_{n-1}$ signed permutations.

\item	If $\Phi_1=i^{(1)}$ then, since $\Phi$ avoids $T$,  
	the symbols $i+1,\dots,n$ are not signed by $1$, and the symbols $1,\dots,i-1$ are not signed by $2$. 
	Hence there are $|E_{n-i}^{r-1}||E_{i-1}^{r-1}((1^{(1)},2^{(1)}))|$ signed permutations, 
	which means by part $(i)$ of Main Theorem that there are 
	${{n-1}\choose{i-1}} (n-i)!(r-1)^{n-i}d_{r-1}(i-1)$ signed permutations.
\end{enumerate}
Therefore $p_n$ staisfies the following relation:
   $$p_n=n(r-1)p_{n-1}+\sum_{i=1}^n {{n-1}\choose{i-1}} (n-i)!(r-1)^{n-i}d_{r-1}(i-1).$$
Besides $p_0=1$, and $p_1=r$, hence the theorem holds.
\qed\end{proof}

\begin{example} ({\em see} \cite[Eq. $48$]{Rs})
Let us denote $a_n=|E_n^2((1^{(1)},2^{(1)}),(2^{(2)},1^{(1)}))|$;  
by symmetric operations and by Theorem \ref{case4},  
$a_n=na_{n-1}+(n-1)!\sum\limits_{j=0}^{n-1} \frac{1}{j!}$ for $n\geq 1$,  
hence $n!<p_n<(n+1)!$ for $n\geq 3$.
\end{example}
%=================================================================
%=================================================================
\begin{corollary}
Let $wc(r)$ be the number of $r$-Wilf classes of a double restriction by $2$-letter 
signed patterns. Then for $r\geq 1$
	$$wc(r)=\left\{ \begin{array}{ll}
			1	&,\ \mbox{if}\  r=1\\
			4	&,\ \mbox{if}\  r=2\\
			6	&,\ \mbox{if}\  r\geq 3
			\end{array} \right..$$
\end{corollary}
\begin{proof}
By Theorem \ref{case1}, Theorem \ref{case2}, Theorem \ref{case3}, and Theorem \ref{case4} 
we get $wc(1)=1$, $wc(2)=4$. The rest follows from definitions, and by the first elements 
of $E_n^3(T)$ where $T$ set of two signed patterns from Table \ref{tab1}.
\qed\end{proof}
%========================================================================
\section*{\centering{\sc 7. Combinatorial identity}}
First of all let us define for $a_1\leq a_2\leq \dots\leq a_{2l}$
	$$U_{a_1,\dots,a_{2l}}^{b_1,\dots,b_l}=\{ (1^{(a_{2i-1})},2^{(a_{2i})}) | b_i=0\} 
					       \cup \{ (2^{(a_{2i-1})},1^{(a_{2i})}) | b_i=1 \},$$
where either $b_i=0$, or $b_i=1$ for $i=1,2,\dots,l$. So by part $(i)$ of the Main Theorem 
we obtain the following corollary.

\begin{corollary}
Let $1\leq a_1\leq \dots\leq a_{2l}\leq r$, and $b_i\in\{0,1\}$ for all $i=1,2,\dots,l$. Then 
 $|E_n^r(U_{a_1,\dots,a_{2l}}^{b_1,\dots,b_l})|=|E_n^r(U_{1,1,2,2,\dots,l,l}^{0,\dots,0})|.$
\end{corollary}

By Theorem \ref{disgood} $|E_n^r(U_{1,1,2,2,\dots,l,l}^{0,\dots,0})|$ is equal to 
$$\sum_{i_1+\dots+i_l\leq n} {n\choose {i_1,\dots,i_l,n-i_1-\dots-i_l}}^2 (n-i_1-\dots-i_l)! (r-l)^{n-i_1-\dots-i_l}.$$
On the other hand, by part $(i)$ of Main Theorem $|E_n^r(U_{1,2,3,4,\dots,2l}^{0,\dots,0})|$ is equal to 
$$\sum_{i_1+\dots+i_l\leq n} {n\choose {i_1,\dots,i_l,n-i_1-\dots-i_l}}^2 (n-i_1-\dots-i_l)! (r-2l)^{n-i_1-\dots-i_l} \prod_{j=1}^l d_2(i_j).$$
where $d_2(m)=\sum\limits j!{m\choose j}^2$. 
Hence by the above Corollary we obtain the following theorem.
\begin{theorem}
\label{id}
Let $r\geq 2l$. For $n\geq 0$ 
  $$
  \begin{array}{l}
  \sum\limits_{i_1+\dots+i_l\leq n} \frac{{n\choose {i_1,\dots,i_l}}^2}{(n-i_1-\dots-i_l)!}(r-l)^{n-i_1-\dots-i_l}=\\
  \ \ \ \ \ \ \ \ \ \ \ \ \ \ \ \ \ \ =\sum\limits_{i_1+\dots+i_l\leq n} \frac{{n\choose {i_1,\dots,i_l}}^2}{(n-i_1-\dots-i_l)!}(r-2l)^{n-i_1-\dots-i_l} \prod_{j=1}^l d_2(i_j).
  \end{array}$$
\end{theorem}

\begin{example}
By Theorem \ref{id} for $l=1$, and $r=3$, we obtain 
   $$\sum\limits_{i=0}^n \frac{2^{n-i}}{i!^2(n-i)!}=\sum\limits_{i=0}^n\sum\limits_{j=0}^i \frac{1}{j!(i-j)!^2(n-i)!}.$$
\end{example}
%**********************************************************************
%**********************************************************************
%=======================================================================

\end{document}